\documentclass{amsart}

\usepackage{amsmath, amsthm}
\usepackage{amssymb, latexsym}
\usepackage{amsfonts}

\newtheorem{Theorem}[equation]{Theorem}
\newtheorem{theorem}[equation]{Theorem}
\newtheorem{Proposition}[equation]{Proposition}
\newtheorem{proposition}[equation]{Proposition}
\newtheorem{lemma}[equation]{Lemma}
\newtheorem{Lemma}[equation]{Lemma}

\newtheorem{corollary}[equation]{Corollary}

\newtheorem{remark}[equation]{Remark}

\newcounter{com}

\newcommand{\beql}[1]{\begin{equation}\label{#1}}
\newcommand{\eeq} {\end{equation}}

\font\Aaa=msam10

\def\qed{\hbox{~~\Aaa\char'003}}

\font\Bbb=msbm10
\def\semi{\hbox{\Bbb o}}

\def\Z{\hbox{\Bbb Z}}

\def\N{\hbox{\Bbb N}}

\numberwithin{equation}{section}

\def\M{{\mathcal M}}

\let\define=\def
\let\redefine=\def

\def\Max{{\rm Max}}

\def\Aut{{\rm Aut}}
\def\Out{{\rm Out}}

\def\G{\Gamma}

\redefine\D{{ \Delta }}

\redefine\O{{ \Omega }}

\define\({ \left( }
\define\){ \right) }
\define\[{ \left[ }
\define\]{ \right] }
\define\<{ \langle }
\define\>{ \rangle }
\let\ljunk=\{
\let\rjunk=\}
\redefine\{{\left\ljunk}
\redefine\}{\right\rjunk}

\redefine\.{ \diamomd }
\redefine\+{\oplus}
\redefine\.{ \cdot}
\redefine\#{\sharp}

\redefine\.{ \cdot }

        \def\lcm{{\rm lcm}}

        \def\remark{\noindent{\bf Remark.~}}
        \def\O{{\rm O}}

        \def\PSL{{\rm PSL}}
        \def\SL{{\rm SL}}

        \def\SU{{\rm SU}}
        
        \def\PSU{{\rm PSU}}

        \def\PSp{{\rm PSp}}
        
        \def\Sp{{\rm Sp}}
        
        \def\GL{{\rm GL}}
      
        \def\Out{{\rm Out}}

        \font\Aaa=msam10

\font\Aaa=msam10

\def\qed{\hbox{~~\Aaa\char'003}}

\font\Bbb=msbm10
\def\semi{\hbox{\Bbb o}}

\def\Z{\hbox{\Bbb Z}}

\def\div{ \kern-.5pt\hbox{\big |} }
\def\ndiv{ {\not\kern-.5pt\hbox{\big |}\,} }
\def\ndivv{ {\not\kern+1.5pt\hbox{$\mid$}\,} }

\def\Aut{{\rm Aut}}

\def\GL{{\rm GL}}

\def\B{^2\kern-.8pt B}
\def\G{^2\kern-.8pt G}
\def\EH{^2\kern-.8pt\hat  E}
\def\E{^2\kern-.8pt E}
\def\D{^3\kern-1pt D}
\def\FF{^2\kern-.8pt F}

\newdimen\refcodesize
\newbox\seriesbox
\def\para#1{\medskip\noindent{\bf #1}~}

\def\proof{\noindent {\bf Proof.~}}

\def\PSL  {{\rm PSL }}

\def\O{{\rm O}}
\def\Po{{\rm P\Omega}}

\def\PSU{{\rm PSU}}
\def\PSp{{\rm PSp}}

\def\Sp{{\rm Sp}}

\def\GL{{\rm GL}}
\def\SU{{\rm SU}}

\def\SL{{\rm SL }}
\def\GL{{\rm GL}}

\def\lcm{{\rm lcm}}

\begin{document}

\title[Invariable generation]
{Invariable generation and the Chebotarev invariant of a finite group}

\thanks{The authors acknowledge partial support from 
  NSF grant DMS~0753640 (W.\hspace{1pt}M.\hspace{1pt}K.),    
ERC Advanced Grants 226135 (A.\hspace{1pt}L.) and  
247034 (A.\hspace{1pt}S.), and ISF grant 754/08 (A.\hspace{1pt}L. and 
A.\hspace{1pt}S.).
The first author is grateful for  the warm hospitality of the Hebrew 
University while this paper was being written.}

       \author{W. M. Kantor}
       \address{University of Oregon,
       Eugene, OR 97403}
       \email{kantor@uoregon.edu}

    \author{A. Lubotzky}
       \address{Institute of Mathematics, Hebrew University, Jerusalem 91904}
       \email{alexlub@math.huji.ac.il}
       
    \author{A. Shalev}
       \address{Institute of Mathematics, Hebrew University, Jerusalem 91904}
       \email{shalev@math.huji.ac.il}

{

}

\begin{abstract}
A subset $S$ of a finite group $G$ {\em invariably generates} $G$ 
if $G=\<s^{g(s)}\mid s\in S  \>$ for each choice of $g(s)\in G, s\in S$.  
We give a tight upper bound on the minimal size of an invariable generating 
set for an arbitrary finite group $G$. In response to a question in 
\cite{KZ} we also bound the size of a randomly chosen set of elements 
of $G$ that is likely to generate $G$ invariably.
Along the way we prove that every finite simple group  is invariably generated by two elements.


 \end{abstract}

\maketitle

\vspace{-8pt}
\centerline{\em Dedicated  to Bob Guralnick in honor  of  his 60th birthday}

\vspace{4pt}

\section{Introduction}

\label{Introduction}

For many years there has been a rapidly growing literature   concerning the generation  of finite groups.  
This has involved the number $d(G)$ of generators of a group $G$, or the 
expected number $E(G)$ of random choices  of elements in order to probably generate $G$, among other group-theoretic invariants.
In this paper we will study further invariants.  
   
Dixon \cite{Di1} began the probabilistic direction for generating (almost) simple groups, and later he also introduced yet another direction based on the 
goal of determining Galois groups
\cite{Di2}.  This has led to the following notions:

\medskip
{\noindent \bf Definition.}
Let $G$ be a finite group.
\begin{itemize}
\item [\rm(a)] A subset $S$ of $G$ {\em invariably generates} $G$ if 
$G=\<s^{g(s)}\mid s\in S  \>$ for each choice of $g(s)\in G, s\in S$ 
\cite{Di2}.
\item [\rm(b)] Let  $d_I(G):=
\min\{  |S|\,\big| \,S \mbox{ invariably generates } G\}$.
\item [\rm(c)] The {\em Chebotarev invariant} $C(G)$ of $G$ is 
the expected value of the random variable $n$  that is minimal subject 
to the requirement that $n$ randomly chosen elements of $G$ invariably generate $G$
\cite{KZ}.
\end{itemize}
\medskip

There have been several papers discussing 
(a) for specific groups (such as finite simple groups) 
\cite{LP,NP,Sh,FG,KZ}, but not for finite groups in general.
Concerning (c), recall Chebotarev's Theorem that provides elements of a 
suitable Galois group $G$,
where the elements are  obtained only  up to conjugacy in $G$; 
the interest in (c) comes from computational group theory, where there is a need to know how long one 
should expect to wait in order to ensure that  choices of representatives from 
the conjugacy classes provided by Chebotarev's Theorem will generate $G$.  
This is discussed more carefully in \cite{Di2,KZ}.

Our main results are the next two theorems, which depend on the  classification of the finite simple groups.

\begin{theorem}
\label{Theorem 1} 
Every finite group $G$ is invariably generated by at most $\log_2|G|$ elements.
\end{theorem}

This bound is best possible:
we show that $d_I(G) = \log_2|G|$ if and only if $G$ is an elementary abelian
$2$-group. It is trivial that $d(G)\le \log_2|G|$
using Lagrange's Theorem. However,
$d_I(G)$ may be much larger than $d(G)$:
Proposition~\ref{powers} states that, {\em for every $r\ge1,$ there is a finite group $G$ such that $d(G)=2$ but $d_I(G)\ge r$.}
Theorem~\ref{composition length} contains a more precise statement 
of Theorem~\ref{Theorem 1} involving the length and structure of a chief series of $G$.

\begin{theorem}
\label{Theorem 2} There exists an absolute constant $c$ such that
$$C(G) \le c|G|^{1/2}(\log|G| )^{1/2}$$
 for all finite groups $G$.
\end{theorem}

This bound is close to best possible: 
it is easy to see that sharply 2-transitive groups provide 
an infinite family of groups $G$ for which $C(G)\sim |G|^{1//2}$ 
(compare \cite[Sec.~4]{KZ}). 
In fact   \cite[Sec.~9]{KZ} asks whether 
$C(G) = O(|G|^{1/2})$ for all finite groups $G$
(which we view as rather likely).

For an arbitrary finite group it is interesting to compare $d_I(G)$ with $d(G)$, and $C(G)$ with $E(G)$.   
The upper bounds for $d_I(G)$ and $d(G)$ are identical, 
although (as stated above) these quantities may be very different.  
On the other hand, $E(G)\le ed(G)+2e\log\log|G|+11 = O(\log|G|)$  ~\cite{Lu}, 
 which is far smaller than the bound in Theorem~\ref{Theorem 2}.

We will  need the following result of independent interest.

\begin{theorem}
\label{Theorem 3}
Every nonabelian finite simple group is invariably generated by $2$ elements. 

\end{theorem}

In fact, for proofs of   Theorems~\ref{Theorem 1}
and \ref{Theorem 2}
 we will need  slightly stronger results on simple groups
involving automorphisms as well (cf. Theorems~\ref{Theorem 3A} 
and \ref{Theorem 3C}).
The same week that we proved these results about simple groups essentially 
the same result as Theorem~\ref{Theorem 3A} with a roughly similar proof was posted in \cite{GM2}.

Dealing with simple groups uses the 
rather large literature of known properties of those groups.
The fact that, for finite simple 
groups $G$,  $d_I(G)$ and $C(G)$ are bounded  by some (unspecified) constant $c$ follows 
for alternating groups  from \cite{LP} 
(cf. \cite{KZ}), and for Lie type groups from 
results announced  in \cite{FG} related to ``Shalev's 
$\epsilon$-Conjecture'', 
which concerns the number of fixed-point-free elements in simple 
permutation groups (cf. Section~\ref{proof of Theorem 2}).

The proof of Theorem~\ref{Theorem 2}  uses  bounds in \cite{CC}
and \cite{FG} on the number of fixed-point-free elements of a transitive 
permutation group, together with a recent bound on the number of maximal 
subgroups of a finite group \cite{LPS}.
We note that an explicit formula for $C(G)$ is given in 
\cite[Proposition~2.7]{KZ},  but we have not been able to use it since 
it appears to be too 
difficult to evaluate its terms for most groups $G$.

 The proofs of Theorems~\ref{Theorem 1},
\ref{Theorem 2} and 
\ref{Theorem 3} are given in Sections~\ref{proof of Theorem 1},
\ref{proof of Theorem 2} and \ref{proof of Theorem 3}, respectively.  Section~\ref{Preliminaries} contains the aforementioned result on the non-relationship of $d(G)$ and $d_I(G)$, as well as a characterization of nilpotent groups as those finite groups all of whose generating sets invariably generate.

\medskip
This paper is dedicated to Bob Guralnick, who 
has made fundamental contributions in the various areas 
involved in 
this and other papers of ours.

\section{Preliminary results and examples}
\label{Preliminaries}

Unless otherwise stated, we assume  that the group $G$ is finite. 
If $X, Y \subseteq G$,  we say that $Y$ is
\emph{similar}  to $X$ if there is a function $f\colon X \rightarrow Y$ such
that $f(X) = Y$ and, for each $x \in X$, $f(x)$ is conjugate in $G$
to $x$. Thus $X$ invariably generates $G$ if and only if 
$\langle Y \rangle = G$ for each $Y \subseteq G$ that is 
similar to $X$.  

Let $\Max(G)$
denote the set of maximal subgroups of $G$. Let ${\M} =\M(G)$ be a set
of representatives of conjugacy classes of maximal subgroups of
$G$. 

If $M\in \Max(G)$, write
$$\widetilde{M} = \bigcup_{g \in G} M^g  \mbox{\ \  and  \   }
v( M ) =  \frac{ |\widetilde{M}|}{|G|}\:.$$
Clearly $\widetilde{M_1} = \widetilde{M_2}$ if the maximal subgroups
$M_1, M_2$ are conjugate in $G$.
Also, $\widetilde{M}$ is
the set of elements of $G$ having at least one fixed point in the
primitive permutation representation of $G$ on the set $G/M$ of 
(left) cosets of $M$ in $G$.

\begin{lemma}
\label{invariable generation criterion}
 A subset $X \subseteq G$ generates $G$ invariably
if and only if $X \not\subseteq \widetilde{M}$ for all $M \in \M$.
\end{lemma} 

\proof If  $X \subseteq \widetilde{M}$ for some $M \in \M$ then
each element of $X$ is conjugate to an element of $M$, and
hence $X$ does not generate $G$ invariably.
Conversely, if $X$ does not generate $G$ invariably, then
there exists a set $Y$ similar to $X$ such that $\langle Y \rangle \ne G$.
Hence (using the finiteness of $G$) there exist $M \in \M$ 
and $g\in G$ 
such that $\langle Y \rangle \subseteq M^g$, and hence
$X \subseteq \widetilde{M}$.
\qed
\medskip

The ``only if'' part of the above lemma also holds for infinite
groups.  Moreover,  the proof shows that
 $X \subseteq G$ generates an arbitrary  group $G$ invariably
 only if $X \not\subseteq \widetilde{H}$ for all $H<G$.
This enables us to show that some
infinite groups are not invariably generated by any set
of elements. For example, there are countable groups $G$ all of whose nontrivial elements are conjugate \cite{HNN}
(and even 2-generated groups  with this property \cite{Os}), so that 
$\widetilde{H}=G$
for every nontrivial subgroup $H$   and hence 
 even $G$ itself does 
not generate $G$ invariably.

However, for finite groups there are no anomalies of this kind,
since $\widetilde{H} \ne G$ for all proper subgroups $H$.
In fact, if  $k(G)$ denotes the number of conjugacy classes of (elements of) 
the finite group $G$, then we have

\begin{lemma}
For any finite group $G$ we have $d_I(G) \le k(G)$. Moreover$,$
$d_I(G)$ is at most the number of conjugacy classes 
of cyclic subgroups of $G$. 
\end{lemma}
\proof
If $H$ is the subgroup of $G$ generated by a set of cyclic subgroups, 
one from each conjugacy class, then 
the union of all conjugates of $H$ is $G$, and hence $H=G$.~\qed
\medskip

For $k \ge 1$, let $P_I(G,k)$ be the probability that $k$ randomly chosen 
elements of $G$ generate $G$ invariably.

\begin{lemma} 
\label{trivial bounds}
$\displaystyle \max_{M\in \M}v(M)^k \le 1-P_I(G,k) \le
 \sum_{M \in \M} v(M)^k $.
\end{lemma}

\proof Let $g_1, \ldots , g_k \in G$ be randomly chosen.
Given $M \in \M$, the probability that $g_i \in \widetilde{M}$ for
all $i$ is $v(M)^k$. Both inequalities now follow
easily from Lemma~\ref{invariable generation criterion}.~\qed
\medskip

We next  characterize nilpotent groups in terms
of invariable generation.

\begin{Proposition}
\label{nilpotent}
A finite group $G$ is nilpotent if and only if every generating set 
of $G$ invariably generates $G$.
\end{Proposition}

\proof
Let $\Phi(G)$ denote the Frattini subgroup of $G$. Then a subset
   of $G$ generates $G$ if and only if its image  in
$G/\Phi(G)$ generates $G/\Phi(G)$.

Suppose $G$ is nilpotent. Then $G/\Phi(G)$ is abelian.
Suppose $X  \subseteq G$ generates $G$, and let
$Y  \subseteq G$ be similar to $X$. Clearly the images of $X$
and $Y$ in the abelian group $G/\Phi(G)$ coincide. Since the image of 
$X$ generates $G/\Phi(G)$, so does the image of $Y$. It follows that
$Y$ generates $G$. We conclude that $X$ invariably generates $G$.

Now suppose $G$ is not nilpotent. We shall construct
a generating set $X$ for $G$ that does not generate $G$ invariably using a theorem of Wielandt \cite[p. ~132]{R}:
if $G/\Phi(G)$ is abelian then $G$ is nilpotent. Then
$G/\Phi(G)$ is not abelian, and hence some maximal
subgroup $M$ of $G$ is not normal   in $G$. Let $g \in G$ with
$M^g \ne M$.  Let $x \in M^g \setminus M$ and $X: = M \cup \{ x \}$.
Then $\langle X \rangle = G$ since $M$ is maximal, so that 
$M \cup \{ x^{g^{-1}} \} = M$ is similar to $X$
and is proper in $G$. This
implies that $X$ does not generate $G$ invariably.
\qed
\medskip

In particular, for nilpotent $G$ we have $d_I(G)=d(G)$.  For simple groups, 
by Theorem~\ref{Theorem 3}  we also have the same equality 
(with both sides 2).  However, our next result shows that,
in general, $d_I(G)$ is not 
bounded above by any function of $d(G)$:

\begin{proposition}
\label{powers}
For every $r\ge1$ there is a finite group $G$ such that $d(G)=2$ but $d_I(G)\ge r$.
\end{proposition}

This group $G$ will be a power $T^k$ of an alternating group $T$.  For this purpose we recall an elementary criterion in \cite[Proposition~6]{KL}:

\begin{proposition}
\label{KL criterion}
Let $G=T^k$ for a nonbelian finite simple group $T$.
Let $S=\{s_1,\dots,s_r\}\subset G$, so that 
$s_i=(t_1^i,\dots,t_k^i), t_j^i\in T$.  Form the matrix
$$
A=\begin{pmatrix} 
        t_1^1&\dots & t_k^1 \\
         &\dots & \\
        t_1^r&\dots & t_k^r \\
     \end{pmatrix}.
     $$
Then $S$ generates $G$ if and only if the following both hold$:$
\begin{itemize}
\item[\rm(a)]
If $1\le j\le k$ then $T=\< t_j^1,\dots,t_j^r  \>;$ and
\item[\rm(b)] The columns of $A$ are in different $\Aut(T)$-orbits for the diagonal action of  $\Aut(T)$ on $T^r$.
\end{itemize}
\end{proposition}

{\noindent \bf Proof of Proposition~\ref{powers}.}
Fix $n$, let $T=A_n$ and let $k=k(n)$ be the 
largest integer such that $d(G)=2$, where $G:=G_n=T^{k}$. 
Then $k\ge n!/8$ \, (\cite[Example~2]{KL}, obtained from Proposition \ref{KL criterion}).

Let $S$ be as in Proposition~\ref{KL criterion}, and assume that $S$ invariably generates $G$.  Then we can arbitrarily conjugate each $t_j^i$ independently and still generate $G$. 
Let ${\bf C}(T)$ denote the set of conjugacy classes of $T$. Project each column $\beta_j$ of $A$ to 
$\bar \beta_j\in {\bf C}(T)^{r}  $.  In view of conditions (a) and (b) in Proposition~\ref{KL criterion}, the $\bar \beta_j$ are in different 
$\Aut(T)$-orbits of the diagonal action on 
${\bf C}(T)^{r}  $.

The number of conjugacy classes in $T$ is at most $c^{\sqrt n}$, 
so  $|{\bf C}(T)|^{r} \le c^{r\sqrt n}$.
The number of projections $\bar \beta_j$ is $k$ (since $1\le j\le k$), 
where $k\ge n!/8$.   Then  $c^{r\sqrt n} \geq  n!/8$ by the Pigeon Hole Principle,
so that $|S|=r \geq C \sqrt n\log n$.~\qed

  

\section{Proof of Theorem~\ref{Theorem 1}}
\label{proof of Theorem 1}

Let $l(G)$ denote the length of a chief series of $G$.
The following is a stronger version of Theorem~\ref{Theorem 1}:

\begin{Theorem}
\label{composition length}
Let $G$ be a finite group having a chief series with $a$ abelian
chief factors and $b$ non-abelian chief factors.  Then
$$d_I(G) \le a + 2b.$$
In particular$,$ $d_I(G) \le 2l(G),$ and if $G$ is solvable
then $d_I(G) \le l(G)$.
\end{Theorem}

\proof
We use induction on $|G|$ (the case $|G|=1$ being
trivial). Suppose $|G| > 1$ and let $N \lhd G$ be a minimal normal
subgroup of $G$. It suffices to show that
$$d_I(G) \le d_I(G/N) + c,$$
where $c=1$ if $N$ is abelian and $c=2$ if $N$ is non-abelian.
In the latter case our proof relies on Theorem~\ref{Theorem 3A}  (proved below).

Let $X \subseteq G$ be a set of size $d_I(G/N)$ whose image in
$G/N$ generates $G/N$ invariably.

Suppose first that $N$ is abelian. Let $x \in N$ be any non-identity
element of $N$. We claim that $Y = X \cup \{ x \}$ 
\emph{invariably generates
$G$}. Indeed, suppose $Z \subseteq G$ is similar to $Y$. 
Then the image of $Z$ in $G/N$ generates $G/N$ (by the assumption on $X$). 
Moreover, $Z$ contains a conjugate $z = x^g$ that is a non-identity element 
of $N$.
Since $G/N$ acts irreducibly on $N$,   $\<Z\>\ge  N$. It follows that $\<Z\>=G$, so $Y$ generates $G$ invariably. Thus
$d_I(G) \le d_I(G/N) + 1$ in this case.

Now suppose $N$ is non-abelian. Then $N = T_1 \times \cdots \times T_k$,
where $k \ge 1$ and the $T_i$ are non-abelian
finite simple groups such that   the conjugation action of $G$ on $N$
induces a transitive action of $G/N$ on the
set $\{ T_1, \ldots , T_k \}$.

The group $A := N_G(T_1)/C_G(T_1)$ is an almost simple group with socle 
$T_1^\star:=T_1C_G(T_1)/C_G(T_1)\cong T_1$.
By Theorem~\ref{Theorem 3A}, there are elements $x_1\in T_1^\star, $  $ x_2 \in A$ such that 
$\langle x_1^{a_1}, x_2^{a_2} \rangle \ge T_1^\star$ for all $a_1, a_2 \in A$.
Let $y_1\in T_1, y_2 \in N_G(T_1)$, be pre-images of $x_1, x_2$,
respectively.  
We claim that $Y: = X \cup \{ y_1, y_2 \}$   {\em invariably generates $G$.}%

To see this, let  $Z$  be a set similar to $Y$, so $Z = X'\cup \{ y_1^{g_1},y_2^{g_2}\}$ 
where $X'$  is similar to $X$ and  
$g_i\in G$ ($i=1,2)$.  We need to show that $Z$ generates $G$.
Let $K = \langle Z \rangle$ and $H = \langle X' \rangle$. 
Since $X$ invariably generates $G$
modulo $N$ we have $HN = G$. Hence $H$ acts transitively
(by conjugation) on $\{  T_1, \ldots , T_k \}$. 

Moreover,  $T_1^{g_1} = T_i$ and $T_1^{g_2} = T_j$ for some $i, j$.
By the transitivity of $H$ there are elements $h_1, h_2 \in H$ 
such that $T_i^{h_1} = T_1$ and $T_j^{h_2} = T_1$.
Then
$g_1h_1, g_2h_2 \in N_G(T_1)$.

Clearly $y_1^{g_1h_1}\in T_1^{g_1h_1}=T_1$ and 
$y_2^{g_2h_2}\in N_G(T_1)^{g_2h_2}= N_G(T_1).$  Then 
$y_1^{g_1h_1}$ and $y_2^{g_2h_2}$ induce automorphisms of $T_1$ by conjugation.
In view of our choice of $x_1$ and $x_2$, 
$\<y_1^{g_1h_1}, y_2^{g_2h_2}\>$ induces all inner automorphisms of $T_1$.
 In particular, the conjugates of the element 
 $y_1^{g_1h_1}\in T_1$ under this group  generate the simple group  $T_1$.  
Thus,  $K\ge \<y_1^{g_1h_1},y_2^{g_2h_2},H  \>\ge T_1$, so that 
$K \ge T_i$
for all $i$ and hence  $G=KN = K$,
as required.  

We see that $d_I(G) \le d_I(G/N) + 2$ in the non-abelian case.
This completes the proof of the first assertion in the theorem.
The last two assertions follow immediately.
\qed
\medskip

{\noindent \bf We can now complete the proof of Theorem~\ref{Theorem 1}.}   Let $G,a,b$ be as above. 
Every abelian chief factor of $G$ has order at least $2$, 
while every non-abelian chief factor has order at least $60$. 
This yields $|G| \ge 2^a 60^b$, so  that 
$$\log_2|G| \ge a + (\log_2{60})b \ge a+2b \ge d_I(G),$$
as required.
Moreover, if $d_I(G)= \log_2|G|$ then we must have $b=0$,
and all chief factors of $G$ have order $2$. Thus $G$ is a 2-group, 
so that $d_I(G)=d(G)= \log_2|G|$ by Proposition \ref{nilpotent}.
Now $d(G) = \log_2 |G|$ easily implies that $G$ is  
an elementary abelian $2$-group.  \qed
\medskip

Note that the bound in Theorem~\ref{composition length} is tight
both for   non-abelian
simple groups and  for
elementary abelian $p$-groups.
\medskip

\section{Proof of Theorem~\ref{Theorem 2}}
\label{proof of Theorem 2}

The main result of this section is the following.

\begin{theorem}
\label{square root theorem}
For any $\epsilon > 0$ there exists $c = c(\epsilon)$
such that $P_I(G,k) \ge 1 - \epsilon$ for any finite group $G$ and any $k \ge c |G|^{1/2}(\log{|G|})^{1/2}$.
\end{theorem}

\proof  
For  $M \le G$ let $M_G=\cap_{g\in G}M^g$ denote the 
{\em core} of $M$ in $G$, the kernel of the permutation action of $G$ on the set of conjugates of $M$.

Divide the set $\M$ of representatives of 
conjugacy  classes of maximal subgroups of $G$
 into three subsets $\M_1, \M_2, \M_3$ as follows.
The set $\M_1$ consists of the subgroups $M \in \M$ such that
the primitive group $G/M_G$ is not of affine type.
The set $\M_2$ consists of the subgroups $M \in \M$ such that
the primitive group $G/M_G$ is of affine type and 
 $|G\colon\! M| \le |G|^{1/2}/(\log{|G|})^{1/2}.$
Finally, $\M_3$ consists of the remaining subgroups in $\M$,
namely the subgroups $M$ such that $G/M_G$ is affine
and  $|G \colon \!M| > |G|^{1/2}/(\log{|G|})^{1/2}.$

By  \cite[Theorem~1.3]{LPS}, for any finite group $G$ we have
$|\Max(G)| \le c_1 |G|^{3/2},$
where $c_1$ is an absolute constant.
In particular, for  $i=1,2,3$,
$$|\M_i| \le |\M| \le  c_1 |G|^{3/2}.$$

Fix $k \ge 1$ and let $g_1, \ldots , g_k\in G$ be randomly chosen
(we will restrict $k$ in later parts of the proof).
By Lemma 2.1,
$$1-P_I(G,k) \le P_1 + P_2 + P_3,$$
where $P_i$ is the probability that $g_1, \ldots , g_k \in \widetilde{M}$
for some $M \in \M_i$ ($i=1,2,3$). It suffices to show that,
for $k$ as in the statement of the theorem, $P_i < \epsilon/3$
for $i = 1,2,3$.

We bound each of the probabilities $P_i$ separately.
By increasing the constant $c$ 
we may assume
that $|G|$ is as large as required in various parts of the proof.

\para{The set $\M_1$.}\  
To bound $P_1$ we use   \cite[Theorem~8.1]{FG}: the proportion of fixed-point-free permutations
in a non-affine primitive group of degree $n$ is at least $c_2/\log{n}$,
for some absolute constant $c_2 > 0$.
This shows that, for $M \in \M_1$,
$$v(M) \le 1 - c_2/\log{|G \colon \!M|} \le 1 -  c_2/\log{|G|}.$$

By Lemma 2.3 and its proof,
$$P_1 \le \sum_{M \in \M_1} v(M)^k \le |\M_1| (1 -  c_2/\log{|G|})^k
\le c_1 |G|^{3/2} (1 -  c_2/\log{|G|})^k.$$
Since $(1-x)^k \le \exp(-kx)$ for $0<x<1$,
 for any $c_3>
\log c_1+3/2$ the right
hand side is bounded above by 
$\exp(c_3 \log{|G|}  - c_2k/{\log{|G|}})$.
If $k >c_4(\log{|G|})^2$ for a suitable absolute  constant $c_4$,  
then the latter expression
tends to zero as $|G| \rightarrow \infty$, and hence so does
$P_1$. In particular we have
$P_1 < \epsilon/3$ for $|G|$ large enough.

\para{The set $\M_2$. }\
We next bound $P_2$. Here our main tool is the theorem  
 that the proportion
of fixed-point-free elements in any transitive permutation
group of degree $n$ is at least $1/n$ \cite{CC}.
This implies that, if  $M \in \M_2$, then
  $$v(M) \le 1-|G:M|^{-1} \le 1-(|G|/\log{|G|})^{-1/2}.$$
Therefore
$$P_2 \! \le  \!\!\sum_{M \in \M_2} \!\!
v(M)^k \le |\M_2| \big ( 1 - (|G|/\log{|G|})^{-1/2} \big)^k
\!   \le c_1 |G|^{3/2} \big(  1 - ( |G|/\log{|G|} )^{-1/2}  \big )^k.$$
As before the right side is bounded above by
$\exp(c_3\log{|G|}  - k \big(|G|/\log{|G|})^{-1/2} ) \big)$
for suitable $c_3>3/2$. 
This in turn tends
to zero as $|G| \rightarrow \infty$ for any 
$k > c_5|G|^{1/2}(\log{|G|})^{1/2}$, for arbitrary $c_5>c_3$.
Therefore $P_2 \rightarrow 0$ for such $k$, and $P_2 < \epsilon/3$
for all sufficiently large $|G|$.

\para{The set $\M_3$. }\
Finally we bound $P_3$. 
If  $M \in \M_3$ then
$G/M_G = V\semi H$, where $V$ is an elementary abelian $p$-group
for some prime $p$, acting regularly on the set of cosets of $M$
in $G$, and $H$ is a point-stabilizer
acting irreducibly on $V$.

Fix a chief series $\{ G_i \}$ of $G$. Fix ${M \in \M_3}$, and let  ${\pi \colon G \rightarrow G/M_G}$ be the canonical projection.
The series $\{ \pi(G_i) \}$
of normal subgroups  of
 $\pi(G) = G/M_G$ descends from $G/M_G = V \semi H$ to $1$.
If  $i$ is minimal  such that $\pi(G_{i+1}) = 1$, then
$\pi(G_i)$ is a minimal normal subgroup of $G/M_G$,
 and hence is $V$, 
the unique minimal normal subgroup of $G/M_G$. In this situation we shall say that
$M$ {\it uses} 
$G_i/G_{i+1}$, in which case  $G_i/G_{i+1} \cong V$. 
(For, since $\pi(G_i)=\pi(G_i)/\pi(G_{i+1})$ is a nontrivial $G$-homomorphic image of $G_i/G_{i+1}$ 
it is isomorphic to $G_i/G_{i+1}$.)
We have seen that every $M \in \M_3$ uses $G_i/G_{i+1}$ for a 
unique $i$. 
Moreover, since $M\in \M_3$,
$$|G_i\!:\!G_{i+1}| = |V|=|G\!:\! M| > (|G|/\log{|G|})^{1/2}.$$

We claim that, {\em if $G$ is sufficiently large$,$ then it has at most two abelian
chief factors used by any maximal subgroups in $\M_3$.} 
Indeed, if there were (at least)
three such chief factors, appearing at places $i>j>l$ in our
chief series, then we would obtain the contradiction
$$|G| \ge |G_i \colon\! G_{i+1}||G_j \colon\! G_{j+1}||G_l \colon\! G_{l+1}| >  
\big((|G|/\log{|G|})^{1/2} \big)^3.$$

Fix an abelian chief factor $V=G_i/G_{i+1}$ of $G$ as above.
Then each $g \in G_i \setminus G_{i+1}$ acts fixed-point-freely
on the cosets of any $M$ that uses $G_i/G_{i+1}$
(since $gM_G \in V \setminus \{ 1 \}$).
For each such $M$ we have
$\widetilde{M} \subseteq G \setminus (G_i \setminus G_{i+1}).$
Since
$$|G \colon \!G_i| \le |G|/|G_i \colon\!G_{i+1}| 
=|G|/|V|\le (|G| \log{|G|} )^{1/2}$$
by the definition of $\M_3$, the proportion of 
\vspace{2pt}
elements $g \in G_i \setminus G_{i+1}$
inside $G$ is at least ${1 \over 2} |G\colon\! G_i|^{-1} \ge {1 \over 2} 
(|G| \log{|G|})^{-1/2}$.
Since the union of $\widetilde M^k$ over all $M$ using
$G_i/G_{i+1}$ is contained in 
$\big (G \setminus (G_i \setminus G_{i+1})\big)^k$,
it follows that
the probability that randomly chosen elements
$g_1, \ldots , g_k$ of $G$ all lie in $\widetilde{M}$
for some such $M$ is at most 
$ (1 -  {1 \over 2} (|G| \log{|G|})^{-1/2} )^k $.
Although there may be many choices for $M$ in $\M_3$, 
there  are at most two choices for the chief factor $G_i/G_{i+1}$.
Thus,  
$$P_3 \le 2 \big(1 -  {1 \over 2} (|G| \log{|G|})^{-1/2} \big)^k
\le 2 \exp \big( -{k \over 2} (|G| \log{|G|})^{-1/2} \big), $$
where the right hand side is 
less than $\epsilon/3$ for $k \ge c (|G| \log{|G|})^{1/2}$
for some $c = c(\epsilon)$.
\smallskip

Our bounds on the three probabilities $P_i$ complete the proof. ~\qed
\medskip

\remark
Recall that the {\em $\epsilon$-conjecture}, posed by the third author of 
this paper,
states that there exists an absolute constant $\epsilon > 0$
such that the proportion of fixed-point-free elements in any 
finite simple transitive permutation group  is at least $\epsilon$.
This amounts to saying that
$v(M) \le 1-\epsilon$
 for any finite simple group $G$ and any   $M\in \Max(G)$.
This conjecture holds for alternating groups \cite{LP} and
for Lie type groups of bounded rank \cite[Secs. 3 and 4]{FG}. Moreover, 
in \cite[Theorem~1.3]{FG} it is announced that the $\epsilon$-conjecture holds
in general, and proofs in some additional cases appear in \cite{FG2}.
When ${M \in \M_1}$ our proof of
Theorem~\ref{square root theorem}  uses  
\cite[Theorem 8.1]{FG}, which in turn relies on the $\epsilon$-conjecture.
However, we now show that Theorem~\ref{Theorem 3C} below easily yields 
a weaker version of 
\cite[Theorem 8.1]{FG} that still suffices for our purpose.

\para{The set $\M_1$ revisited.}\  
Namely, we claim that there exists $c_2>0$ such that 
$$v(M) \le 1 - c_2 (\log |G|)^{-2} |G|^{-1/3},$$
{\em where $G$ is any non-affine primitive permutation group 
and $M$ is a point-stabilizer.}
For, if $s_1, s_2$ generate $G$ invariably, and if  
${M\in\Max(G),}$ then $\widetilde {M} \cap s_i^G = \emptyset$
for   $i = 1$ or $2$, in which case
$v(M) \le 1 - |s_i^G|/|G|.$ 
Then $v(M) \le 1 - {1 \over 2}|G|^{-1/3}$
for  each sufficiently large finite simple group $G$ and 
each such $M$,
by Theorem~\ref{Theorem 3C}. This implies that,  for   all  finite simple groups $G$ and all
  $M\in \Max(G)$, we have  $v(M) \le 1 - c_3|G|^{-1/3}$
  for some constant
$c_3 > 0$.

Consequently, if $G$ is an almost simple group with socle $T$ then,
since $|\Out(T)| \le c_4 \log|T|$ 
(cf. \cite[Sec.~2.5]{GLS}), we easily obtain
$$v(M) \le 1 - c_5 (\log |G|)^{-1} |G|^{-1/3}$$
for all  $M\in \Max(G)$ not containing $T$, for some $c_5 > 0$.
Our claim follows by combining this inequality with the reduction to 
almost simple groups given in the proof of  \cite[Theorem~8.1]{FG}.

Thus, if $M \in \M_1$, then the above claim yields
$$P_1 \le \sum_{M \in \M_1} v(M)^k \le 
c_1 \log|G|  (1 - c_2 (\log |G|)^{-2} |G|^{-1/3})^k.$$
The right hand side tends to zero when $k \ge c_6 (\log|G|)^3|G|^{1/3}$; 
but for the proof of Theorem~\ref{square root theorem}  we 
can assume the stronger inequality 
$k \ge c_7 |G|^{1/2} (\log|G|)^{1/2}$.
Consequently $P_1\to 0$, as required. 

\para{Completion of proof of Theorem \ref{Theorem 2}}. Apply
Theorem \ref{square root theorem} with $\epsilon = 1/2$
and let $c = c(1/2)$. Let $k=\lceil c|G|^{1/2} (\log|G|)^{1/2}\rceil$. Then $k$
randomly chosen elements of $G$ invariably generate $G$
with probability at least $1/2$. This implies that
$$C(G) \le  2k \le(2c+1)|G|^{1/2} (\log|G|)^{1/2}. \ \ \ \qed$$ 

\begin{corollary}
\begin{itemize}
\item[(a)] If $G$ is a finite group without abelian composition factors$,$  
then $C(G) = O((\log{|G|})^2)$.

\item[(b)] If $G$ is an almost simple group$,$ then
$C(G) = O(\log{|G|} \log \log |G|)$.
\end{itemize}
\end{corollary}

\proof 
We have already seen   (a) in  our first treatment of the 
non-affine case ($M\in \M_1$) of   Theorem~\ref{square root theorem}.

To prove (b) we first note that, for some $c > 0$ and 
all $M \in \M$, we have $v(M) \le 1 - c/\log|G|$.
Indeed, if $M$ has trivial core then this follows from
\cite[Theorem 8.1]{FG} (and hence 
 from  the correctness of the $\epsilon$-conjecture stated above). Otherwise, $M$ contains the simple socle
 $T$ of $G$, and $|G/T| \le  |\Out(T)| \le c_4 \log|T| \le c_4 \log|G|$
as noted above. In this situation, if $g \in G$ acts fixed-point-freely on the cosets of 
$M$ in $G$, so do all the elements of $gT$, so that $v(M) \le  1 - c_4^{-1}/\log|G|$. 

By
\cite[Theorem~1.3]{GLT},  $|\M| \le c_1 (\log{|G|})^3$ when $G$ is almost
simple. This yields
$$\sum_{M \in \M} v(M)^k \le c_1 (\log{|G|})^3 (1- c/\log|G|)^k
\le c_1 (\log{|G|})^3 \exp({- ck/\log|G|} ).$$
The right hand side tends to zero as $|G| \rightarrow \infty$ when
$k \ge c_2 \log|G| \log \log |G|$. This proves part (b).
\qed
\medskip

We observe that {\em the bound in} (b) {\em is almost 
best possible, up to the $\log \log{|G|}$ factor.}
To show this we use the following example \cite[p.~115]{FG}.
Fix any prime $p$. Let $G = \PSL(2,p^b).b$, the extension of the simple group by the group $B$
of $b$ field 
automorphisms, where $b$ is a prime not dividing $p(p^2-1)$.
Let $G$ act on the cosets of the maximal subgroup $N_G(B)$ of $G$. 
Then all fixed-point-free elements are contained in the socle
of $G$, so their proportion is less than $1/b$.
Therefore $v(M) \ge 1- 1/b$.

Hence, by Lemma~\ref{trivial bounds},    
$P_I(G,k) \le 1- (1-1/b)^k$,
 so that for sufficiently large  $b$  we obtain
$$P_I(G,k) \le 1 - (1-c_1/\log|G|)^k \le 1 - \exp(-c_2k/\log|G|),$$ 
where $c_1, c_2$ are suitable constants.
Thus $P_I(G,k) \le 1/2$ for all $k \le c_3 \log|G|$,
where $c_3 > 0$ is an absolute constant.
The probability that it takes at least $k+1$ random choices of elements
to invariably generate $G$ is $1-P_I(G,k)$.
By the definition of the expectancy $C(G)$ we have 
$C(G) \ge (k+1)(1-P_I(G,k))$. If $k=[c_3 \log|G|]$ then $1-P_I(G,k) \ge 1/2$ 
and $k+1 \ge c_3 \log|G|$.
This yields
$C(G) \ge (k+1) (1 / 2) \ge (c_3/2)\log|G|.$

\section{Simple groups}
\label{proof of Theorem 3}

We will prove the following slightly stronger version of Theorem~\ref{Theorem 3}:

\begin{Theorem}
\label{Theorem 3A}
Let $G$ be a finite simple group.  
\begin{itemize}
\item[\rm(a)]
 If $G$ is not one of the groups $\Po^+(8,q),$ $q=2$ or $3,$ then there are two  elements
$s_1,s_2\in G$ such that 
$G=\<s_1^{g_1}, s_2^{g_2} \>$ for each choice of 
$g_i\in \Aut(G).$
 
\item[\rm(b)]
 If $G$ is  $\Po^+(8,q),$ $q=2$ or $3,$  and if $G\le G^\star\le \Aut(G),$  then there are  elements
$s_1\in G,s_2\in G^\star$ such that 
$G\le \<s_1^{g_1}, s_2^{g_2} \>$ for each choice of 
$g_i\in G^\star.$
 
\end{itemize}
\end{Theorem}
  
 Of course, Theorem~\ref{Theorem 3} is just (a) using inner automorphisms.
 This theorem is also obtained in \cite[Theorem~7.1]{GM2}, 
 along with the fact that
  $\Po^+(8,2) $ is an actual exception.  
 
 We begin with the easiest case:

\begin{Lemma}
\label{alternating}
{\rm Theorem~\ref{Theorem 3A}} holds for 
each alternating group $A_n,$ $n\ge5$.
  
\end{Lemma}
\proof
If $n\ne6$ then $\Aut(A_n)=S_n$.
 For even  $n>6$  use the  product  of a disjoint
$2$-cycle and   $(n-2)$-cycle, and the product of a disjoint
$p$-cycle and   $(n-p)$-cycle for a prime $p\le n-3$ not dividing $n$; it is easy to check that such a prime exists.  These two elements generate a group $H$ 
that is readily seen to be transitive and even primitive.
Since $H$ contains a $p$-cycle, $H=A_n$ by  
a classical result of Jordan \cite[Theorem~13.9]{Wie}.

If $n$ is odd then an $n$-cycle and a $p$-cycle can be used in the same manner, 
for an odd prime $p\le n-3$ not dividing $n$.

Finally,   $A_6$ is generated by any elements of order 4 and 5.   \qed

\medskip
For groups of Lie type we will use the knowledge of all maximal overgroups  $M$  of a carefully chosen  semisimple element  $t_1$.   
Then,
by Lemma~\ref{invariable generation criterion},
we only need to  choose an
 $\Aut(G)$-conjugacy class of elements that does not  
 meet
 the union of the corresponding  sets $\widetilde  M$.
Our arguments differ from those in \cite{GM2}
primarily due to that paper using \cite{GM} whereas we rely more on  the earlier paper \cite{MSW}.

\begin{Lemma}
\label{classical}
{\rm Theorem~\ref{Theorem 3A}} holds for 
each classical simple group  
other than $\Po^+(8,q)$.
\end{Lemma}

{
\font\sevenroman=cmr8
\font\seventemp=cmsy8
\font\sevenital=cmmi8
\textfont0=\sevenroman
\textfont2=\seventemp
\textfont1=\sevenital

\begin{table}[t] 
  \caption{Classical groups}
  \label{classical torus}
  \vspace{-18pt}
$$\begin{array}{|l|l|l|l|l|}    
\hline
 \mbox{\sevenroman  quasisimple } G\!\!& |t_1|&
 t_1 ~\mbox{\sevenroman  on }V&  |t_2|& t_2 ~\mbox{\sevenroman  on }V\\
 \hline 
\SL(n,\!q) & (q^n-1)/(q-1)& n&(q^{n-1}-1)/(q-1)
&(n-1)\oplus1
\\\ \  n \text{ \sevenroman odd}&&&&
\\ \hline
\SL(n,\!q) 
&(q^{n-1}\!-\!1)/(q\!-\!1)\!\!
&(n-1)\oplus1
& (q^n-1)/(q-1)& n
\\\ \  n\ge4\text{ \sevenroman even}&&&&\\ 
\hline
\Sp(2m,q)& q^m+1&2m& \lcm(q^{m-1}+1,q+1)
&(2m-2)\perp 2\!\!
\\\ \  m\ge2 &&&&
\\ \hline
\Omega(2m+1,q)& (q^m+1)/2
& 2m^-\perp 1  &  (q^m-1)/2& (m\oplus m)\perp 1 
\\\ \  q \text{ \sevenroman odd}&&&&
\\ \hline
\Omega^+(4k,q)
 &( q^{n'-1}+1)/\delta_1 \raisebox{2.3ex} {~} 
 &(n-2)^-\!\!\perp\! 2^-\!\!\!
&\lcm(q^{n'-2}\!+\!1,q^{2}\!+\!1)/\delta_2\!\!
& (n-4)^-\perp 4^-\!\!
\\ \,\, n=2n'=4k\!&&&&
 \\ \hline
\Omega^+(4k+2,q)\!
 &( q^{n'-1}+1)/\delta_1 \raisebox{2.3ex} {~}  
 &(n-2){}^-\!\!\perp 2^-\!\!& 
(q^{n'}-1)/\delta_2
& n'\oplus n'
\\ \,\, 2n'\!=4k+2\!&&&&
 \\ \hline
\Omega^-(4k,q)& 
(q^{n'}+1)/\delta_1 \raisebox{2.3ex} {~} 
&n^-  & (q^{{n'}-1}-1)/\delta_2
&(n-2)^+\perp 2^- \!\!
\\ \,\, n=2n'\!=4k&&&&
\\ \hline
\Omega^-(4k+2,q)\!&
 (q^{2k+1}+1)/\delta_1 \raisebox{2.3ex} {~} 
 &(4k+2)^-  & 
(q^{2k}+1)/\delta_2
&4k^-\perp 2^+  
\\ \hline
\SU(2m,q)
&q^{2m-1}+1&(2m-1)\perp 1
&(q^{2m}-1)/(q+1)&2m
\\ \hline
\SU(2m+1,q) \!& (q^{n}+1)/(q+1)&n&q^{n-1}-1&n-1\perp 1
\\ \hline
\end{array}
$$
\end{table}
} 
\proof
We will  consider the corresponding quasisimple  linear group $G$,
using    semisimple elements $t_1$ and   $t_2$   in Table~\ref{classical torus}  
 that  decompose the space as indicated in the table.
 (Here $\delta_i$ is  1 or 2,
$n$ is the dimension of the underlying vector  space $V$, and $n'=n/2$. 
If an entry involves $\lcm(q^i+1,q^j+1)$ for some $i,j$, then 
$t_2$ induces irreducible elements of order $q^i+1$  or  $q^j+1$ on the indicated subspaces of dimension $2i$ or $2j$.)

In each case, $t_1$ is the element called  ``$s$''   in \cite[Theorem~1.1]{MSW}; if there is a $1-$ or $2-$space indicated then it is 
centralized. 
 For each group $G$, all 
maximal overgroups of $t_1$ 
are listed in  \cite[Theorem~1.1]{MSW}.
Until the end of the proof we will exclude the case $G=\Sp(4,q)$. 
Then all automorphisms of $G$ act on $V$, preserving the underlying geometry \cite[Sec.~2.5]{GLS}.  It follows that
all $\Aut(G)$-conjugates of $t_i$ act on $V$ as $t_i$ does (for $i=1,2$).
We always use conjugates of $t_1$ and $t_2$ that have no assumed relationship to one another, so if the two elements studied generate $G$ then they invariably generate $G$.

If $G$ is not $\SL(2,q)$, $\Sp(4,q)$  or $\Sp(8,2)$,  then  $t_1$ and $t_2$ invariably generate 
$G$ by \cite[Theorem 1.1]{MSW}:
all of the exceptions in that theorem do not arise here due to the behavior of {\em both} $t_1$ and $t_2$ on $V$.  If $G=\Sp(8,2)$  then we   replace $t_2$ by another element, as follows.  Let $f\in G$ have order 5 and centralize a nondegenerate $4-$space.  
Then $C_G(f)=\<f\>\times \Sp(4,2)$.  Let $c=(1,2)(3,4,5,6)\in 
S_6\cong \Sp(4,2)<C_G(f)$.  Then $c\notin S_5\cong \O^-(4,2)$, and hence
$fc$ is not in an overgroup  $\O^-(8,2)$ of $t_1$. Since its order
implies that $fc$  is also not in any of the other maximal overgroups of $t_1$
\cite[Theorem~1.1]{MSW}, it follows that $t_1$ and $fc$ invariably generate  $G$.  

Case $\SL(2,q)$.  When $q$ is $4, 5$ or 9, see Lemma~\ref{alternating}.  When $q=7$, elements of order 7 and 4 invariably generate $G$.  For all other $q\ge4$, the same $t_1$ and $t_2$  as indicated in the table (but with $t_1$ acting irreducibly on each $1-$space) invariably generate $G$ by \cite[Ch.~XII]{Di}.

Case   $\Sp(4,q)$.  We  may assume that $q\ge4 $ since 
 $\Sp(4,2)$  is not simple and $\PSp(4,3)\cong \PSU(4,2)$. 
We again use $t_1$ and $t_2$ as in the table, such that
$t_2$ induces  an element of order $q+1$ inside the $\Sp(2,q)$
produced by each factor in the decomposition $4=2\perp 2$.     Once again 
$t_1$ and $t_2$  invariably generate $G$ by  \cite[Theorem~1.1]{MSW}.~\qed

\medskip

We note that classical groups were considered in \cite[Section~10]{NP} from a probabilistic point of view:  a large number of pairs of elements was described that invariably generate various classical groups. The group $\GL(n,q)$ was also handled in \cite{Sh} for large $n$.  All groups of Lie type also were dealt with probabilistically, at least for bounded rank, in \cite[Theorem~5.3]{FG}. 

\begin{Lemma}
\label{O+8}
{\rm Theorem~\ref{Theorem 3A}} holds for 
$\Po^+(8,q)$.

\end{Lemma}

\proof
Once again we will consider the corresponding linear group  $G=\Omega^+(8,q)$, using the properties of $\Aut(G/Z(G))$ contained in \cite[Sec.~2.5]{GLS}.  We have 
$G/Z(G) \break\le G^\star\le \Aut(G/Z(G))$.

(a) Suppose first that $q>3$.   We will  use the same $\<t_1\>$ as above (mod $Z(G)$), of order  $ (q^3+1)/(2,q-1)$.
It acts on our space as $8^+=6^-\perp 2^-$,
centralizing the $2-$space.  
 
We also use an element $t_3 \in G$ of order $(q^3-1)/(2,q-1)$.  Here $t_3 $ 
decomposes our space as  $8^+=(3\oplus 3)\perp (1\oplus 1)$
using totally singular $3-$ and $1-$spaces, inducing  isometries  of order 
$q-1$ on the subspace
$1\oplus 1$
and of order 
$q^3-1$ on the subspace $3\oplus 3$, and hence acting irreducibly on the indicated 
$3-$spaces.  Then  $t_3 $ fixes exactly two singular $1-$spaces, and two 
totally singular $4-$spaces in each $G$-orbit of such 4-spaces
(each of the latter fixed subspaces has the form $3\perp 1$).
If $\tau $ is any   automorphism of $G/Z(G)$, then 
$t_3 ^\tau$ has the same properties.  In particular, neither $t_3 $ nor
$t_3 ^\tau$ fixes any anisotropic $1-$ or $2-$space
for any $\tau\in \Aut(G/Z(G))$. 
(N.\:B.--This requires that $q>3$:  if $q=3$ then the analogous element $t_3 $ induces $-1$ on the $2^+-$space $1\oplus 1$ and hence fixes all of its $1-$spaces.)
However,  by   \cite[Theorem 1.1]{MSW} each maximal subgroup 
of $G/Z(G)$ that 
contains $t_1$ (mod $Z(G))$ either fixes such a $1-$ or $2-$space or its
image under a triality automorphism behaves that way.  Hence, there is no maximal subgroup containing  $t_1 $ and $t_3$ mod $Z(G)$, and we have invariably generated $G/Z(G)$.

\smallskip
(b) From now on $q\le 3$.  
First consider the case where $G^\star$ acts (projectively)  on $V$ (this includes the situation in Theorem~\ref{Theorem 3}). 
 We use   elements  $t_3$ and $t_4$ of $G/Z(G)$ of order $(q^4-1)/(4, q^4-1)$
arising from a decomposition $8^+=4^-\perp 2^-\perp 2^+$
and  from a decomposition $8^+=4\oplus 4$ into totally singular 
$4-$spaces  (the corresponding cyclic groups $\<t_i\>$ are conjugate under 
$\Aut(G/Z(G))$ but not under $G^\star$).
The Sylow $5$-subgroups of $\<t_3\>$ and $\<t_4\>$   behave differently on the vector space, and $\<t_4\>$ 
is an element
of order $(q^4-1)/(4, q^4-1)$ that
acts fixed-point-freely on $V$. Hence, by 
\cite{Kl},   
$\<t_3,t_4\>$ is contained in no proper subgroup of $G$,
so that $\<t_3,t_4\> = G$.

Finally, suppose that  $G^\star$  does not have any
  $\Aut(G/Z(G))$-conjugate  that acts on   $V$.
  Here we return to the original setting of  the theorem, now letting $G$ denote the simple group $\Po^+(8,q)$.
  Since Out$(G)\cong S_4$ or $S_3$,  we may assume that $G^\star$ contains a  triality outer automorphism.
Consequently,  there is a subgroup
$\Z_3\times\SL(3,q)$  of  $G^\star$  
that contains an element $t_5$ of order 
$3(q^2 + q +1)$ such that $\tau=t_5^{q^2+q+1}$ is  a triality automorphism   and $t_5^3$ acts projectively on $V$ as 
$8^+=(3\oplus 3)\perp(1\oplus 1)$.   

By \cite[Theorem 1.1]{MSW},  $\<t_1,t_5\>\cap G  \ge 
\<t_1,t_5^3\>$ 
 is either $G$, $\Omega(7,q)$ 
 or lies in $A_9<\Omega^+(8,2)$.
 Since $\<t_1,t_5\>\cap G $  is invariant
  under $\tau$, only the first of these can occur (for example, 
 $\<t_1,t_5\>\cap G $  cannot be  $A_9$ or $\PSL(2,8)<A_9$).
 Thus, $t_1$ and  $t_5$ invariably generate $G \<\tau\>$.~\qed

\medskip

{\noindent \bf Completion of proof.} 
In
\cite[Tables~6~and~9]{GM} there are lists of carefully chosen  cyclic subgroups of exceptional and sporadic simple groups, as well as all of the maximal 
overgroups $M$ of those subgroups.  It is straightforward to use those tables to handle these final cases of Theorem~\ref{Theorem 3A}. 
This amounts to exhibiting an element order for $G$ not appearing in any of the listed subgroups $M$.
We provide some details for the exceptional groups.  
 Table~\ref{exceptional torus} reproduces part of \cite[Table~6]{GM}. 
 Here $T_1$ is a cyclic maximal torus and $M$ runs through the isomorphism types of  maximal overgroups of   $T_1$.
(Notation:  $\epsilon =\pm1$, $\Phi_n =\Phi_n (q)$ is the $n$th cyclotomic polynomial
evaluated at~$q$, $\Phi_8'=\Phi_8'(q)=q^2+\sqrt{2}q+1$,
$\Phi_{12}'=\Phi_{12}'(q)=q^2+\sqrt{3}q+1$  and  
$\Phi_{24}'=\Phi_{24}'(q)=q^4+\sqrt{2}q^3+q^2+\sqrt{2}q+1$.)
 In each case,   the order of $t_2$ guarantees that it is not contained in any of the listed maximal overgroups $M$ (there are also other choices for $t_2$).
Hence, a generator  of $T_1$ together with  $t_2$ behave as required in the theorem. \qed  

\medskip

In Section~\ref{proof of Theorem 2}   we needed a bit more information than in the preceding theorem for an 
alternative proof of Theorem~\ref{square root theorem} 
and hence of Theorem 1.2: 
\begin{Theorem}
\label{Theorem 3B}
\label{Theorem 3C}
For all sufficiently large $G$ in 
{\rm Theorem~\ref{Theorem 3A}}$,$ the elements $s_i$   can be chosen so that 
$|s_i^G|>|G|^{2/3} /2$ for $i=1,2$.  
\end{Theorem}

\proof
This is a straightforward matter of examining each 
part of the proof of Theorem~\ref{Theorem 3A}.
In each case we need to check that $|C_G(s_i)|< 2|G|^{1/3}\,$
for $i=1,2$ and   all sufficiently large $|G|$.

For alternating groups, when $n$ is even each
 of the groups $C_G(s)$ 
is the direct product of two cyclic groups, and hence has order satisfying the required bound.  When $n$ is odd the same holds if we 
replace the $p$-cycle  by  the product of a disjoint $p$-cycle  and an  $(n-p)$-cycle 
(a power of which is a $p$-cycle). 

In Lemma~\ref{classical} --
excluding $\SL(2,q)$~-- we have 
$|C_G(T_1) | \sim q^r$ and $|C_G(t_2) | \sim q^r$,  where $r$ is the rank of the corresponding algebraic group.  (For example, 
for $\SL(n,q) $ we have 
$|C_G(T_1) |  = (q^n-1)/(q-1)$
or
$q^{n-1}-1$,
for $\Sp(2m,q)$ we have 
$|C_G(t_2) |  \le (q^{m-1}+1)  (q+1)$, and for 
$\Omega^+(4k+2,q)$ we have $|C_G(T_1) | \le (q^{2k}+1)(q+1)$.)
A straightforward calculation using $|G|$  verifies that these bounds are  small enough for our purposes.  When $G=\SL(2,q)$ we have 
$|C_G(T_1) | =q+1$, so that $|s_i^G|>|G|^{2/3} /2$ and  a denominator larger than 1 is essential. 
\qed

{
\font\sevenroman=cmr8
\font\seventemp=cmsy8
\font\sevenital=cmmi8
\textfont0=\sevenroman
\textfont2=\seventemp
\textfont1=\sevenital

\begin{table}[t] 
  \caption{Exceptional groups}
  \label{exceptional torus}
  \vspace{-18pt}
$$\begin{array}{|l|l|l|l|l|}    
\hline
 G& |T_1|&  M\ge T_1& \mbox{\sevenroman further max.}
 & |t_2|\\
 \hline
 ^2B_2(q^2)& \Phi_8'&  N_G(T_1)& - 
 \raisebox{2ex} {~}\raisebox{-1.2ex} {~} &\Phi_8'(-q)
 \\ 
 \ \ q^2\ge8  & &   &   
 &  \\\hline
 ^2G_2(q^2)\!\!& \Phi_{12}'&  N_G(T_1)& 
  -\raisebox{2ex} {~}\raisebox{-1.2ex} {~}& 
 \Phi_{12}'(-q)
 \\ 
 \ \ q^2\ge27  & &   &   
 &  \\\hline 
 G_2(q),\: 3|q+\epsilon\!\!
 & q^2+\epsilon q+1
 \raisebox{2.2ex} {~}
&  \SL ^\epsilon(3,q).2& \PSL(2,13) 
 & q^2-\epsilon q+1\\ 
 & &   &  \  (q=4)
 &  \\\hline
 G_2(q),\: 3|q& q^2+q+1&  \SL (3,q).2& \PSL(2,13) 
 & q^2- q+1\\
 & &   &  \   (q=3)
 &  \\\hline
   ^3D_4(q)& \Phi_{12}&  N_G(T_1)& -
   \raisebox{2.2ex} {~}\raisebox{-1.2ex} {~}
   & (q^3 \!+\! 1)(q \!- \!1)/(2,q\!-\!1)\!\! \\ \hline
 ^2F_4(q^2)& \Phi_{24}'&  N_G(T_1)& -
  \raisebox{2.2ex} {~}\raisebox{-1.2ex} {~}
 & \Phi_{24}'(-q) 
 \\ 
 \ \ q^2\ge8  & &   &   
 &  \\\hline
 F_4(q)
 &  \Phi_{12}& ^3D_4(q).3 & \PSL(4,3).2_2,
 & q^4+1 \\
& &   &  ^2F_4(2) ~  (q\!=\!2),\!\!
 & 
 \\
 & &   & \PSL(4,3).2_2  \,  (q\!=\!2)\!\!
 & 
  \\\hline
E_6(q)& \Phi_9/(3,q-1)&  \SL(3,q^3).3& -
       & 
       (q\! + \!1)(q^5 \!-\! 1)/(6,q\!-\!1)\!\!\\
       \hline
^2\hspace{-.5pt}E_6(q)& \Phi_{18}/(3,q+1)&  \SU(3,q^3).3& -
   \raisebox{2.2ex} {~}\raisebox{-1.2ex} {~}& 
     (q \!-\! 1)(q^5\! +\! 1)/(6,q\!+\!1)\!\!\\
    \hline
E_7(q)& \Phi_2\Phi_{18}/(2,q\!-\!1)\!\!& ^2\!E_6(q)_{sc}.D_{q+1} \!\!
\!  & -
       & \Phi_7/(2,q-1)\!\!    
       \\ \hline
E_8(q)& \Phi_{30} & N_G(T_1)& -
       &\Phi_{24}  \\ \hline
\end{array}
$$
\end{table}
}

\para{Random generation.} 
We conclude with remarks concerning the
 random generation of finite
simple groups. All finite simple groups $G$
are generated by two randomly chosen elements with probability
tending to 1 as $|G| \rightarrow \infty$   \cite{Di1,KL,LS}. We claim that this does
not hold for invariable generation: 
{\em the probability that 
two~--~or any bounded number of~--~random elements of a finite simple group $G$ invariably generate $G$ is bounded away from $1$.}
To show this we need the following result that is implicit in
\cite{FG}.

\begin{lemma} There exists an absolute constant $\epsilon > 0$
such that any finite simple group $G$ has a maximal subgroup $M$
for which $v(M) \ge \epsilon$.
\end{lemma}

\proof This is trivial for alternating groups $A_n$, where we
take $M$ to be a point-stabilizer in the natural action, so
$v(M) \sim 1-e^{-1}$. For groups $G$ of Lie type of bounded
rank over a field with $q$ elements we may assume $q$ is large, and then   the result follows with $M$ a maximal subgroup containing a
maximal torus (see the discussion in
\cite[start of Sec.~4]{FG}). For classical groups of large rank
the result follows from    \cite[Theorem~1.7]{FG}. Sporadic simple groups
 satisfy the conclusion trivially.
\qed
\medskip

This lemma can be considered as a kind of weak analogue
of the   $\epsilon$-conjecture (stated above)  but in the opposite direction.     

We can now deduce

\begin{corollary} There is an absolute constant $\epsilon > 0$ such that
$P_I(G,k) \le 1 - \epsilon^k$ for all 
finite simple groups $G$ and positive integers $k$.\end{corollary}

\proof This follows by combining the above lemma with
Lemma~\ref{trivial bounds}.
\qed
\medskip

In \cite[p.~114]{FG} it is announced that, for any $\epsilon > 0$,
there is $c = c(\epsilon)$ such that 
 $P_I(G,k) \ge 1-\epsilon$
 whenever  $G$
is a finite simple group of Lie type and $k \ge c$.
The case of bounded rank is proved in \cite[Theorem 4.4]{FG}, and
a similar result for alternating groups was proved earlier in
\cite{LP}. 

Using these results it follows that,
for any function $f\colon \!\N \rightarrow \N$ such that ${f(n) \rightarrow
\infty}$ as $n \rightarrow \infty$ (even if arbitrarily slowly),
 we have
$P_I(G,f(|G|)) \rightarrow 1$ for finite simple groups $G$
whose orders tend  to infinity.


\end{document}